\PassOptionsToPackage{hyphens}{url}
\documentclass[final]{siamart171218} 

\usepackage{amsmath,amssymb}
\usepackage{color}
\usepackage{upref}
\definecolor{hotpink}{rgb}{0.9,0,0.5}
\hypersetup{colorlinks,urlcolor=blue,citecolor=hotpink,linkcolor=blue}
\usepackage{graphicx}
\usepackage{enumitem}
\usepackage{booktabs}

\usepackage{color}
\renewtheorem{algorithm}[theorem]{Algorithm}

\newcommand{\RR}{\mathbb{R}}

\title{Accurate Computation of the Log-Sum-Exp and Softmax Functions%
       \thanks{Version of September 8, 2019.
       \textbf{Funding}: This work was supported by Engineering and
        Physical Sciences Research Council grant EP/P020720/1, The MathWorks, 
        and the Royal Society.
        The opinions and views
        expressed in this publication are those of the authors, and not
        necessarily those of the funding bodies.}}

\author{Pierre Blanchard%
        \thanks{Department of Mathematics,
                University of Manchester,
                Manchester, M13 9PL, UK
                (\texttt{pierre.blanchard00@gmail.com}).
               }
        \and 
        Desmond J. Higham%
        \thanks{%
           School of Mathematics,
           University of Edinburgh,
           Edinburgh, EH9 3FD, UK
           (\texttt{d.j.higham@ed.ac.uk}).%
                 }
        \and
        Nicholas J. Higham%
        \thanks{%
                Department of Mathematics,
                University of Manchester,
                Manchester, M13 9PL, UK
                (\texttt{nick.higham@manchester.ac.uk}).%
               }
}

\definecolor{pred}{rgb}{0.5,0,0}
\definecolor{pgreen}{rgb}{0,0.5,0}
\definecolor{pblue}{rgb}{0,0,0.6}

\def\new#1{#1}

\def\py{polynomial}

\def\Alg{Algorithm}

\def\pert{perturbation}
\def\fe{forward error}
\def\be{backward error}
\def\rea{rounding error analysis}
\def\reas{rounding error analyses}

\def\eu{\ensuremath{\mathrm{e}}}

\def\xmax{x_{\max}}
\def\xmin{x_{\min}}
\def\rmax{r_{\max}}
\def\rmins{r_{\min}^{(s)}}
\def\rmin{r_{\min}}

\def\d{\delta }
\def\ds{\Delta s}
\def\st{\widetilde{s}}

\def\nrm{normwise}

\def\cn{condition number}
\def\R{\mathbb{R}}
\def\eps{\epsilon}

\def\fp{floating-point}
\def\fpa{floating-point arithmetic}

\def\norm#1{\|#1\|}

\def\normi#1{\|#1\|_1}
\def\normo#1{\|#1\|_{\infty}}

\DeclareMathOperator{\fl}{\operatorname{f\kern.2ptl}} 
\def\op{\mathbin{\mathrm{op}}}

\DeclareMathOperator{\cond}{cond}
\DeclareMathOperator{\condo}{cond_{\infty}}
\DeclareMathOperator{\logonep}{log1p}

\def\C{\mathbb{C}}
\def\nbyn{n \times n}

\def\shat{\widehat{s}}
\def\what{\widehat{w}}

\def\yhat{\widehat{y}}

\def\ghat{\widehat{g}}

\def\us{u} 

\def\lse{log-sum-exp}
\def\sm{softmax}

\def\t#1{\texttt{#1}}
\def\alg{algorithm}
\def\th{\theta}

 \makeatletter
 \def\mymatrix#1{\null\,\vcenter{\normalbaselines\m@th
     \ialign{\hfil$##$\hfil&&\quad\hfil$##$\hfil\crcr
       \mathstrut\crcr\noalign{\kern-\baselineskip}
       #1\crcr\mathstrut\crcr\noalign{\kern-\baselineskip}}}\,}
 \makeatother

\newcounter{mylineno}
\makeatletter
\let\oldtabcr\@tabcr

\def\mynewline{\refstepcounter{mylineno}%
               \llap{\footnotesize\arabic{mylineno}\hspace{5pt}}%
              }

\gdef\@tabcr{\@stopline \@ifstar{\penalty%
           \@M \@xtabcr}\@xtabcr\mynewline}

\newenvironment{code}{%
                        \mathcode`\:="603A  
                        \def\colon{\mathchar"303A}
                        \setcounter{mylineno}{0}
                        \par
                        \upshape
                        \begin{list} 
                           {} {\leftmargin = 1cm}
                        \item[]
                        \begin{tabbing}

                           \hspace*{.3in} \= \hspace*{.3in} \=
                           \hspace*{.3in} \= \hspace*{.3in} \= \kill
                           \mynewline
                       }{\end{tabbing}\end{list}}
\makeatother

\mathchardef\Gamma="7100 \mathchardef\Delta="7101
\mathchardef\Theta="7102 \mathchardef\Lambda="7103
\mathchardef\Xi="7104 \mathchardef\Pi="7105 \mathchardef\Sigma="7106
\mathchardef\Upsilon="7107 \mathchardef\Phi="7108
\mathchardef\Psi="7109 \mathchardef\Omega="710A

\mathcode`@="8000 
{\catcode`\@=\active\gdef@{\mkern1mu}}

\begin{document}
\maketitle

\begin{abstract}
Evaluating the
log-sum-exp function or the softmax function
is a key step in many modern data science algorithms, notably in inference
and classification. 
Because of the exponentials that these functions contain, the evaluation is
prone to overflow and underflow, especially in low precision arithmetic.
Software implementations commonly use alternative formulas
that avoid overflow and reduce the chance of harmful underflow,
employing a shift or another rewriting.
Although mathematically equivalent, these variants behave differently in
floating-point arithmetic.
We give rounding error analyses of different evaluation algorithms
and interpret the error bounds using condition numbers for the functions.
We conclude, based on the analysis and numerical experiments, that the
shifted formulas are of similar accuracy to the unshifted ones and that
the shifted softmax formula is typically more accurate
than a division-free variant.
\end{abstract}



\section{Introduction}

In many applications, especially in a wide range of machine learning
classifiers such as multinomial linear regression and naive Bayes
classifiers 
\cite{cgp19}, \cite{murp12}, \cite{wiba98}, 
one needs to compute an expression of the form
\begin{equation}\label{yold}
        y = f(x) = \log \sum_{i=1}^n \eu^{x_i},
\end{equation}
where $x = [x_1,x_2,\dots, x_n]^T\in\C^n$.  
The function $f:\C^n\to\C$ is often referred to
as \lse\ or LSE\@. Its gradient $g:\C^n\to\C^n$,
given by 
\begin{equation}\label{eq.sm}
  g_j(x) = \frac{\partial}{\partial x_j} f(x) 
         = \frac{\eu^{x_j}}{\sum_{i=1}^n \eu^{x_i}}, \quad
            j=1\colon n,
\end{equation}
is called \sm\ 
and is also a key function in classification algorithms
\cite[p.~355]{efha16}, \cite[p.~78]{gbc16}, \cite{hh19}. 
It is often the case that both \lse\ and \sm\ are required simultaneously. 

The most obvious danger in evaluating \eqref{yold} and \eqref{eq.sm}
is overflow.
We are interested in IEEE arithmetic in the precisions 
half (fp16),
single (fp32),
and double (fp64) \cite{ieee08},
as well as the bfloat16 half precision format~\cite{inte18}.
Table~\ref{table.IEEE} shows the key parameters of interest for these
precisions: the unit roundoff $u$,
the largest finite number $\rmax$, and the smallest positive 
normalized and subnormal \fp\ numbers.
If some $x_i$ exceeds the relevant $\log \rmax$ value in 
Table~\ref{table.exp-max} then overflow will occur.
Clearly, overflow is possible even for quite modestly sized $x$, especially
for half and single precision.

Underflow is also possible.
For example, for $n = 1$, if 
$x_1$ is a finite \fp\ number with 
$x_1 < \log\rmins$ then\footnote{$\log 0 = -\infty$ is the value recommended
by the IEEE standard \cite[p.~43]{ieee08}.}
$\fl(f(x_1)) = \fl(\log(\fl(\eu^{x_1}))) = \fl(\log 0) = -\infty$, 
whereas $f(x_1) = x_1$.
\new{For $n>1$, underflow in the exponential evaluations is a problem when
the sum of the terms that underflow is significant compared with the
sum of the other terms;
otherwise underflows are harmless.}
As well as avoiding harmful underflow, it is desirable to 
avoid generating subnormal numbers, 
which incur a performance penalty if handled in software\footnote{%
\url{https://devblogs.nvidia.com/cuda-pro-tip-flush-denormals-confidence/},
\url{https://en.wikipedia.org/wiki/Denormal_number}.};
see \cite{high:ASNA2} or \cite{mbdj18} for details of subnormal numbers.

A way to avoid overflow, and to attempt to avoid underflow and subnormal
numbers, in evaluating \lse\ is to rewrite
\begin{align*}
          y &= \log \sum_{i=1}^n \eu^{x_i}
             = \log \sum_{i=1}^n \eu^a \eu^{x_i - a}
             = \log  \left(\eu^a\sum_{i=1}^n \eu^{x_i - a}\right).
\end{align*}
If $a\in\R$ then
\begin{equation}\label{ynew}
       y = a + \log \sum_{i=1}^n \eu^{x_i - a}.
\end{equation}
Here and throughout, $\log$ denotes the principal logarithm: the logarithm
whose imaginary part lies in $(-\pi,\pi]$.
Equation~\eqref{ynew} is not, in general, 
true for $a\in\C$ \cite[Lem.~2.5]{aphi14}.
The \sm\ can be expressed in a related form (for any $a$):
\begin{equation}\label{gshift}
      g_j = \frac{\eu^{x_j-a}}{ \sum_{i=1}^n\eu^{x_i-a} },
             \quad j=1\colon n.
\end{equation}
This shifting, typically with $a = \max_i x_i$,
is a well known way to attempt to avoid overflow and
underflow in the evaluation of $f$ and $g$, described in many places,
including on 
Wikipedia\footnote{\url{https://en.wikipedia.org/wiki/LogSumExp}}, in 
blog posts\footnote{For example,
\url{https://hips.seas.harvard.edu/blog/2013/01/09/computing-log-sum-exp/},
\url{http://bayesjumping.net/log-sum-exp-trick/}, and
\url{https://jblevins.org/log/log-sum-exp}. And similarly for the \sm: \url{https://timvieira.github.io/blog/post/2014/02/11/exp-normalize-trick/}.},
and even in a YouTube video\footnote{\url{https://youtu.be/-RVM21Voo7Q}}.
The functions 
\t{logsumexp} in SciPy 1.3.1 \cite{scipy}
and \t{LogSumExp} in R \cite{rlang} 
both implement \eqref{ynew} with $a = \max_i x_i$.
The function \t{softmax} in the 
MATLAB Deep Learning Toolbox (R2019a)
\cite{matlab-DLT} uses \eqref{gshift}
with $a = \max_i x_i$.


\begin{table}[t]
\caption{Parameters for bfloat16 and IEEE fp16, fp32, and fp64 arithmetics, 
         to three significant figures:
         unit roundoff $u$,
         smallest positive (subnormal) number $\rmins$, 
         smallest positive normalized number $\rmin$, and
         largest finite number $\rmax$.
         In Intel's bfloat16 specification,  subnormal numbers are not
         supported, so $\rmins = \rmin$ {\upshape \cite{inte18}}.}

\label{table.IEEE}
\footnotesize


\begin{center}
\begin{tabular}{lllll}\toprule
& \multicolumn{1}{c}{$u$} & \multicolumn{1}{c}{$\rmins$} 
& \multicolumn{1}{c}{$\rmin$} & \multicolumn{1}{c}{$\rmax$}\\\midrule
bfloat16 & $3.91\times 10^{-3}$  &  $9.18\times 10^{-41}$ 
         & $1.18\times 10^{-38}$ & $3.39\times 10^{38}$\\
fp16  & $4.88\times 10^{-4}$ &
$5.96\times 10^{-8}$ & $6.10\times 10^{-5}$ & 
$6.55 \times 10^4$\\
fp32     & $5.96\times 10^{-8}$ &
$1.40\times 10^{-45}$ & $1.18\times 10^{-38}$ & $3.40\times 10^{38}$\\
fp64  & $1.11\times 10^{-16}$ &
$4.94\times 10^{-324}$ & $2.22\times 10^{-308}$ & $1.80\times 10^{308}$\\
\bottomrule
\end{tabular}
\end{center}

\end{table}

\begin{table}
\caption{Logarithms of key parameters in Table~\ref{table.IEEE}, to three
         significant figures.}
\label{table.exp-max}

\label{table.IEEE2}
\footnotesize


\begin{center}
\begin{tabular}{cccc}\toprule
                 & $\log\rmins$ & $\log\rmin$ & $\log\rmax$ \\\midrule
bfloat16 & $-92.2$ & $-87.3$ & $88.7$\\
fp16 & $-16.6 $& $-9.70$ &  $11.0$\\
fp32 & $-103$  & $-87.3$ &  $88.7$\\
fp64 & $-744$  & $-708$  &  $710$   \\
\bottomrule
\end{tabular}
\end{center}
 
 
 
 


\end{table}

An alternative to \eqref{gshift},
which removes the denominator of \eqref{eq.sm} by subtracting \lse\ from
the argument of $\exp$ in the numerator,
is
\begin{equation}\label{galt}
       g_j = \exp\left(x_j - \log\sum_{i=1}^n\eu^{x_i}\right).
\end{equation}
The conciseness of this division-free formula makes it attractive 
for implementing \sm\ when a \lse\ function is 
available.
This formula is used in 
the SciPy 1.3.1 function \t{softmax},
in a MATLAB toolbox \cite{matlab-prml}
associated with the book \cite{bish06},
and in the internal function \t{softmax} in the 
MATLAB Statistics and Machine Learning Toolbox (R2019a) \cite{matlab-SMLT};
in each case the \lse\ term is computed by 
\eqref{ynew} with $a = \max_i x_i$.
The formula \eqref{galt} can also be found in codes posted 
in online communities such as Stack Exchange.


The accuracy properties of the formulas above are not clear.  In
particular, when $a = \xmax < 0$, $y$ in \eqref{ynew} is computed as a sum
of two terms of opposite sign, so there could potentially be damaging
subtractive cancellation.

In this work we analyze the unshifted and shifted formulas
and \eqref{galt} in order to determine which choices of formulas give the best
accuracy and reliability.
In particular, we carry out a \rea\ of \alg s for the evaluation 
and relate the error bounds to the conditioning of $f$ and $g$.
We show that the shifted formulas have broadly similar error bounds to
the unshifted ones, and so are entirely appropriate for practical use.
We find, however, that 
the alternative softmax formula \eqref{galt} has a less favorable
error bound than the shifted formula and tends to produce larger errors in
practice.


We begin, in the next section, by investigating the conditioning 
of the \lse\ and \sm\ functions.
In section~\ref{sec.basic} we give detailed \reas\ of the basic formulas.
In section~\ref{sec.shifted} we analyze the shifted formulas
and \eqref{galt}  and compare their error bounds with those for 
unshifted formulas. Numerical experiments are given in section~\ref{sec.numexp}
to test the accuracy of the evaluations and also to examine how the sum of the 
computed \sm\ vector entries compares with the exact value $1$.
Conclusions are given in section~\ref{sec.conc}.

From this point on, we assume that the $x_i$ are real and 
we write
\begin{equation}\label{xmaxmin}
    \xmax = \max_i x_i, \qquad
    \xmin = \min_i x_i.
\end{equation}
We will use the standard model of \fpa\ \cite[sec.~2.2]{high:ASNA2} 
\begin{equation}\label{eq:fpmodel}
  \fl(a \op b) = (a \op b)(1 + \d), \quad |\d| \le u, \quad
         \op\in\{+,-,\times,/\}.
\end{equation}

\section{Condition number}\label{sec.cond}

Before considering \alg s for computing \lse\ and \sm\
we investigate the conditioning of these functions, 
that is, the sensitivity of $f(x)$ and $g(x)$ 
in \eqref{yold} and \eqref{eq.sm} to small \pert s in $x$.

We define the \cn\ of $f$ in the usual way (see, e.g.,
\cite[chap.~3]{high:FM}), by
\begin{equation*}
      \cond(f,x) := \lim_{\eps\to0} \sup_{\norm{e} \le \eps \norm{x}}
                    \frac{|f(x+e) - f(x)|}{\eps|f(x)|}.
   \label{cond-relf}
\end{equation*}
This definition implies that 
\begin{equation}\label{f-pert}
  \frac{|f(x+e) - f(x)|}{|f(x)|}
  \le \cond(f,x) \frac{\norm{e}}{\norm{x}} + o(\norm{e}),
\end{equation}
so that $\cond(f,x)$ measures the worst-case relative change in $f$
corresponding to a small relative change in $x$.
It is easy to show that for the $\infty$-norm,
\begin{equation}\label{cform}
 \condo(f,x) = \frac{\normi{\nabla f(x)} \normo{x}}{|f(x)|}
            = \frac{\normo{x}}{|f(x)|}
            = \frac{\max_i |x_i|}%
                 {|\log \sum_i \eu^{x_i}|},
\end{equation}
since $\normi{\nabla f(x)} = 1$ by \eqref{eq.sm}.

We identify two extreme cases.  First, the \cn\ is infinite for
$x_i \equiv -\log n$, because $f(x) = 0$.  Hence when $x_i \approx - \log n$
for all $i$ the \cn\ must be large.  Second, if $\max_i x_i = \max_i |x_i|$
then 
$|f(x)| \ge \max_i |x_i|$ by \eqref{ybounds} below, so
$\condo(f,x) \le 1$ and the problem is perfectly conditioned.

A forward stable \alg\ for computing \lse\  is one for which the
relative error of the computed result is bounded by 
$p(n)\cond(f,x)u$, for some low degree \py\ $p$.
Ideally, we would like the \alg\ that we use to 
be forward stable.  
To see whether it is reasonable to expect forward
stability, consider the case $n = 1$.
Then $f(x) = \log \eu^x = x$, so $\cond(f,x) = 1$: the problem is 
perfectly conditioned.
When we compute $f$ using standard library functions we can expect 
to obtain relative errors in the computed exponential and logarithm
bounded by $u$ \cite{dlm07}, \cite{mull16}, \cite[Chap.~10]{mbdj18}, that is,
\begin{equation}\label{yhat}
  \yhat = \fl(f(x)) = \log ( \eu^x (1+\d_1) ) (1+\d_2),
  \quad |\d_1|, |\d_2| \le u.
\end{equation}
The term $1+\d_2$ just causes a small relative \pert\ of the output,
so we have
\begin{align*}
 \yhat \approx  \log ( \eu^x (1+\d_1) ) &= x + \log(1+\d_1) = x + \d_1 + O(\d_1^2).
\end{align*}
Hence, since $y = x$,
\begin{equation}\label{y1b}
  \frac{|y - \yhat|}{|y|} \lesssim \frac{u}{|x|} + O(u^2).
\end{equation}
This relative error bound is much larger than $u$ for $|x| \ll 1$, even
though the problem is perfectly conditioned.
So it is not reasonable to expect an \alg\ to be
unconditionally forward stable in \fpa.
For this trivial computation, \be\ and \fe\ are the same, so 
we also conclude that we cannot expect to obtain an \alg\ that is 
unconditionally backward stable.



 
The \sm\ function has \cn\
\begin{equation*}
      \cond(g,x) := \lim_{\eps\to0} \sup_{\norm{e} \le \eps \norm{x}}
                    \frac{\norm{g(x+e) - g(x)}}{\eps\norm{g(x)}},
   \label{cond-relg}
\end{equation*}
which is given explicitly by 
$$
  \cond(g,x) = \frac{\norm{G(x)}\,\norm{x}}{\norm{g(x)}}.
$$
Here, the $\nbyn$ matrix $G(x) = (\partial g_i/\partial x_j)$ is the
Jacobian of $g$ and $\| \cdot \|$ denotes any vector norm and the
corresponding subordinate matrix norm.
Now
$$
  \frac{\partial g_i}{\partial x_j}
  =
\begin{cases}
     \dfrac{-\eu^{x_i} \eu^{x_j} }{ \biggl(\displaystyle\sum_{k=1}^n \eu^{x_k} \biggr)^2
     },  & i \ne j, \\[30pt]
     \dfrac{ \eu^{x_i} \displaystyle\sum_{k=1}^n \eu^{x_k}  - \eu^{2x_i} }%
          { \biggl(\displaystyle\sum_{k=1}^n \eu^{x_k} \biggr)^2 }, & i = j.
\end{cases}
$$
We have, for each $i$,
\begin{align*}
  \sum_{j=1}^n \left| \frac{\partial g_i}{\partial x_j} \right|
   &=
  \displaystyle\frac{2\eu^{x_i} \displaystyle\sum_{j=1\atop j\ne i}^n \eu^{x_j} }
       { \biggl(\displaystyle\sum_{k=1}^n \eu^{x_k} \biggr)^2}
    \le 1,
\end{align*}
that is, $\normo{G(x)} \le 1$. Hence 
\begin{equation*}
  \condo(g,x) \le \frac{ \normo{x} }{ \normo{g(x)} } \le n\normo{x},
\end{equation*}
because $\normo{g} \ge n^{-1} \normi{g} = n^{-1}$.
We note in passing that $G$ is the Hessian of $f$ and can be shown to be
symmetric positive semidefinite for all $x$~\cite[p.~74]{bova04}.

We also note that shifting, as in \eqref{ynew} and \eqref{gshift}, does not
change the functions so does not change their \cn s; likewise for
\eqref{galt}.  These reformulations may, of course, affect the accuracy of
the \fp\ evaluation.


\section{Basic algorithms and error analysis}\label{sec.basic}

\Alg~\ref{alg.basic} gives a naive implementation of \eqref{yold} and
\eqref{eq.sm}.

\begin{algorithm}\label{alg.basic}
Given $x\in\R^n$, 
this \alg\ computes $f(x) = \log \sum_{i=1}^n \eu^{x_i}$ and
the gradient $g(x) = \nabla f(x)$.
\end{algorithm}
\begin{code}
$s = 0$\\
for \= i = 1:n\\
\> $w_i = \exp(x_i)$\\
\> $s = s + w_i$\\
end\\
$f = \log(s)$\\
for \= i = 1:n\\
\> $g_i = w_i/s$\\
end
\end{code}

What can be said about the accuracy of this \alg\ when it is implemented
in \fpa?
To answer this question we carry out a \rea.
Throughout this section, we assume that there is no overflow or underflow.

First,
we consider the error in evaluating the sum of nonnegative terms
$$
     s = \sum_{i=1}^n \eu^{x_i} \equiv \sum_{i=1}^n w_i.
$$
Evaluating $w_i = \eu^{x_i}$ 
yields a computed result satisfying
\begin{equation}\label{wihat}
     \what_i = \eu^{x_i}(1+\d_1),
\end{equation}
where, as noted in Section~\ref{sec.cond},
we can expect the relative error from the exponential evaluation to satisfy
$|\d_1| \le u$.
Therefore
$$
     |\what_i - w_i| \le w_i u.
$$
Write the (exact) sum of computed quantities as
$$
          \st = \sum_{i=1}^n \what_i.
$$
The rounding error analysis in \cite{high93s}, 
\cite[sec~4.2]{high:ASNA2} shows that the computed sum $\shat$ satisfies
$$
    |\st-\shat| \le \us \sum_{i=1}^{n-1} |t_i| + O(\us^2),
$$
where $t_i = \sum_{j=1}^{i+1} \what_j$, so that,
since $\what_i \ge 0$,
$$
    |\st-\shat| \le \us(n-1)( \what_1 + \what_2 )   
                  + \us\sum_{i=3}^n (n+1-i) \what_i + O(\us^2).
$$
Writing $s - \shat = s - \st + \st - \shat$, we obtain
\begin{align}
  \left|s - \shat\right| 
  &\le \sum_{i=1}^n |\what_i-w_i| + \left|\st - \shat \right| 
      \nonumber\\
  &\le   u \sum_{i=1}^n w_i + \us \sum_{i=1}^n (n+1-i) \what_i
           + O(u^2)\nonumber\\
  &= \sum_{i=1}^n ( n+2-i ) w_i + O(u^2), \label{errb}
\end{align}
since $\what_i = w_i + O(u)$.
Hence
\begin{equation}\label{ds}
    \shat = s + \ds, \quad |\ds| \le (n+1)u s + O(u^2).
\end{equation}
Then the computed \lse\ is
\begin{align}\label{fhat}
 \yhat &= \fl( \log \shat ) = \log(\shat)(1+\eps), \quad |\eps| \le u, \nonumber\\
       &= \log(s+\ds) (1+\eps)\nonumber\\
       &= \left( \log s + \frac{\ds}{s} + O(u^2) \right) (1+\eps)\nonumber\\
       &= y (1+\eps) + \frac{\ds}{s} + O(u^2).
\end{align}
Using \eqref{ds} we obtain
$$
 |y-\yhat| \le u |y| + (n+1)u + O(u^2),
$$
which gives the following result.

\begin{theorem}[Basic \lse\ algorithm]\label{thm.lse-basic}
In the absence of overflow and underflow, 
the computed \lse\ $\yhat$ from \Alg~\ref{alg.basic} satisfies
\begin{equation}\label{yerr2}
 \left|\frac{y-\yhat}{y}\right| 
      \le 
        \left(1 + \frac{n+1}{|y|}\right)u + O(u^2).
\end{equation}
\end{theorem}

Comparing this bound with $\cond(f,x)u$ in \eqref{cform}
we see that it is larger by the factor
$(|y| + n + 1)/\normo{x}$.
But $|y| \le \normo{x} + \log n$ by \eqref{ybounds} below, so this factor is
bounded by $1 + (n+1+\log n)/\normo{x}$.
Hence we have forward stability as long as $\normo{x} \gtrsim 1$, but for
$\normo{x} \ll 1$ the bound does not guarantee forward stability.
This is consistent with the bound \eqref{y1b} for the case $n = 1$.

Turning to the evaluation of the \sm\ function $g$ from 
its definition \eqref{eq.sm}, by \eqref{wihat} we have
$$
   \ghat_j = \frac{ \eu^{x_j} (1+\d_1) }{ \shat } (1+\d_2), 
    \quad |\d_2| \le u, 
$$
where $\d_2$ accounts for the division, and so by \eqref{ds}, 
$$
   \ghat_j = \frac{ \eu^{x_j} }{ s(1 + \eta) } (1+\d_1) (1+\d_2),
             \quad |\eta| \le (n+1)u + O(u^2).
$$
Therefore
\begin{align*}
   \ghat_j &=  g_j (1 + \th), \quad
                |\th| \le (n + 3)u + O(u^2).
\end{align*}
This bound guarantees a relative error of order
at most $nu$ in every component of $g$.
We weaken the bound into a \nrm\ bound for the next theorem.

\begin{theorem}[Basic softmax algorithm] \label{thm.g}
In the absence of overflow and underflow,
the computed \sm\ $\ghat$ from \Alg~\ref{alg.basic} satisfies
\begin{equation}\label{gerr}
  \frac{ \normo{g-\ghat} }{ \normo{g} }
  \le (n + 3)u + O(u^2).
\end{equation}
\end{theorem}



While the error bounds of Theorem~\ref{thm.lse-basic} and \ref{thm.g}
have a very satisfactory form, they provide no useful information when 
$n \gtrsim 1/u$,
 and for fp16 this happens for $n$ as small as $2048$.
We note, however, that the $n$ terms,
which come from the summation, are pessimistic.
It is shown by Higham and Mary~\cite[Thm.~3.1]{hima18a}
that, under a probabilistic model of 
rounding errors, 
$n$ in the error bound for summation  can be replaced by 
a small constant multiple of $\sqrt{n}$ with high probability, and the same 
holds for the bounds of Theorem~\ref{thm.lse-basic} and \ref{thm.g}.


Next, consider the 
alternative formula \eqref{galt}, which we rewrite here:
\begin{equation}\label{galt2}
       g_j = \exp\left(x_j - \log\sum_{i=1}^n\eu^{x_i}\right) = \exp(x_j - y).
\end{equation}
With $y = f(x)$ evaluated in \fpa\ by \Alg~\ref{alg.basic}, we obtain
\begin{align}
 \ghat_j &= (1+\d) \exp\bigl[ (x_j - \yhat) (1+\eps)\bigr], \quad |\d|,|\eps|\le u,  
    \nonumber\\
         &= (1+\d) \exp\bigl[ (x_j - y + (y - \yhat)) (1+\eps) \bigr]  
\label{eq:step2}         
         \\
         &= (1+\d) g_j \exp[ (x_j - y)\eps + (y - \yhat) (1 + \eps) \bigr]
           \nonumber\\
         &= (1+\d) g_j\bigl[(1 + (x_j - y)\eps + (y - \yhat) (1 + \eps) + O(u^2) \bigr]
           \nonumber \\
         &= (1+\th) g_j,
           \label{eq:stepn}
\end{align}
where,
using Theorem~\ref{thm.lse-basic},
$$
      |\th| \le (|y| + |x_j-y| + n + 2 ) u + O(u^2).
$$
We summarize this result as follows.

\begin{theorem}[Alternative softmax algorithm]\label{thm.gshift2}
In the absence of overflow and underflow, the computed $\ghat$ from
\eqref{galt2} with the \lse\ computed by \Alg~\ref{alg.basic} satisfies
\begin{equation}\label{gerr2}
 \frac{ \normo{g-\ghat} }{ \normo{g} } 
 \le \Bigl(|y| + \max_j |x_j-y| + n + 2 \Bigr)u + O(u^2).
\end{equation}
\end{theorem}

From \eqref{ybounds} and \eqref{yxj} below, using the notation
\eqref{xmaxmin}, we have 
$$
 |y| + \max_j |x_j-y|  
      \le |\xmax| + |\xmax - \xmin| + 2 \log n.
$$
Hence \eqref{gerr2} is less favorable than \eqref{gerr} when 
$\xmax - \xmin  \gg n$ or $| \xmax | \gg n$.
The analysis therefore suggests that \eqref{eq.sm} should be preferred to
\eqref{galt}.

To give an intuitive explanation for the potential inaccuracy in (\ref{galt2}), 
we refer to the steps leading to 
(\ref{eq:stepn}).
A large absolute error in the argument of the final exp may lead to a large relative error in the result. This effect can be traced back to the appearance of $x_j -y$ in 
(\ref{eq:step2}).


\section{Algorithms with shifting}\label{sec.shifted}

Now we consider the use of shifts in the \lse\ and \sm\ evaluations
in order to avoid overflow and reduce the chance of \new{harmful} underflow.

Recall the definition \eqref{xmaxmin} of $\xmax$ and $\xmin$.
Overflow in the exponential evaluations in \eqref{ynew} is certainly
avoided if we take
$
   a = \xmax 
$, 
as we then have $x_i-a\le 0$ and hence $0 \le \eu^{x_i-a}\le 1$ 
for all $i$. 
We can rewrite \eqref{ynew} as
\begin{equation}\label{ynew2}
    y = \xmax + \log \Biggl(1 + \sum_{i=1 \atop i\ne k}^n 
        \eu^{x_i - \xmax}\Biggr),
\end{equation}
where $x_k = \xmax$.
From this expression we see that 
\begin{equation}\label{ybounds}
      \xmax \le y \le \xmax + \log n.
\end{equation}
It follows that when $\xmax \ge 0$,
the sum ``$\xmax + \log(\cdot)$'' that produces $y$ cannot suffer cancellation.

Note that for $n = 1$, \eqref{ynew2} trivially 
provides the exact result $y = \xmax$, in contrast to the basic formula \eqref{yold}.


For later use, we note that \eqref{ybounds} implies that, for any $j$,
\begin{equation}\label{yxj}
 |y - x_j| \le |\xmax - x_j| + \log n 
            \le |\xmax - \xmin| + \log n.
\end{equation}

The $\log$ term in \eqref{ynew2} has the form $\log(1+z)$,
where $z \ge 0$.
If $z$ is very small then $1+z$ will round to $1$ 
and the logarithm will evaluate as zero, even though
$\log(1+z) \approx z \ne 0$.
To avoid this loss of information we will
use the function $\logonep(z) = \log(1+z)$ provided in,
for example, C, MATLAB, and 
Numpy.
These functions guarantee an accurate result
for small $z$ 
(which can be achieved with a simple formula based on $\log$
\cite{hp82}, \cite[Prob.~1.5]{high:ASNA2}). 

These considerations lead to \Alg~\ref{alg.log-sum-exp}.


\begin{algorithm}[\lse\ and \sm\ with shift]\label{alg.log-sum-exp}
This \alg\ computes $f(x) = \log \sum_{i=1}^n \eu^{x_i}$ and
the gradient $g(x) = \nabla f(x)$ for $x\in\R^n$.
\end{algorithm}

\begin{code}
$\left[a,k\right] = \max_i x_i$ \% $a = x_k = \max_ix_i$ \\
$s = 0$\\
for \= $i = 1:n$\\
\> $w_i = \exp(x_i-a)$\\
\> if $i \ne k$, $s = s + w_i$, end\\
end\\
$f = a +  \logonep(s)$\\
for \= i = 1:n\\
\> $g_i = w_i/(1 + s)$\\
end
\end{code}

Note that while it is important to avoid forming $1+s$ for the
$f$-evaluation, for $g$ we can safely form $1 + s$ because if $s$ is small
it has little influence on $g$.

\Alg~\ref{alg.log-sum-exp} avoids overflow.
If underflow occurs in the exponential then it is in a term in the sum
added to $1$ in \eqref{ynew2}, so that term is negligible and the underflow
is harmless.
\new{Note, in particular, that if $x_i \approx x < \log \rmins$ for all $i$
then whereas \Alg~\ref{alg.basic} returns $f = -\infty$,
\Alg~\ref{alg.log-sum-exp} suffers no underflow and returns $f \gtrsim \xmax$.}

The main question is how shifting affects the accuracy of the evaluations.
We give a \rea\ to assess this question.
The analysis is a generalization of that in the previous section for the
unshifted \alg.


We 
first examine
the error in evaluating the sum of nonnegative terms
\begin{equation}\label{sum}
    s = \sum_{i=1 \atop i\ne k}^n \eu^{x_i - a}
      =: \sum_{i=1 \atop i\ne k}^n w_i.
\end{equation}
Evaluating $w_i = \eu^{x_i-a}$ 
yields a computed result satisfying
$$
     \what_i = \eu^{ (x_i-a)(1+\d_1) }(1+\d_2), 
      \quad \mbox{$|\d_1|\le u$, $|\d_2|\le u$}.
$$
Therefore
$$
     \what_i = \eu^{ x_i-a } \eu^{(x_i-a)\d_1} (1+\d_2)
             = \eu^{ x_i-a } \bigl(1 + (x_i-a)\d_1 + O\bigl(\d_1^2\bigr)\bigr) (1+\d_2),
$$
and hence 
$$
     |\what_i - w_i| \le ((1+a - x_i)u + O(u^2)) w_i.
$$
Assuming for notational simplicity that $k = n$, we can write the (exact)
sum of computed quantities as
$$
          \st = \sum_{i=1}^{n-1} \what_i.
$$
The rounding error analysis in \cite{high93s}, 
\cite[sec~4.2]{high:ASNA2} shows that the computed sum $\shat$ satisfies
$$
    |\st-\shat| \le \us \sum_{i=1}^{n-2} |t_i| + O(\us^2),
$$
where $t_i = \sum_{j=1}^{i+1} \what_j$, so that,
since $\what_i \ge 0$,
$$
    |\st-\shat| \le 
                   \us\sum_{i=1}^{n-1} (n-i) \what_i + O(\us^2).
$$
Hence
\begin{align}
  \left|s - \shat \right| 
  &\le \sum_{i=1}^{n-1} |\what_i-w_i| + \left|\shat - \st \right| 
      \nonumber\\
  &\le u \sum_{i=1}^{n-1} (1+a-x_i) w_i + \us \sum_{i=1}^{n-1} (n-i)
    \what_i + O(u^2)\nonumber\\
  &= \sum_{i=1}^{n-1} \left(\us\left(n-i\right) +
    u\left(1+a-x_i\right)\right) w_i + O(u^2),  \label{errb2}
\end{align}
since $\what_i = w_i + O(u)$.
%
%
Hence
\begin{equation}\label{werrb}
  \left|\frac{ \shat  - s }{s}\right| \le (n + \xmax - \xmin) u + O(u^2),
\end{equation}
which guarantees an accurate computed sum as long as 
$n + \xmax - \xmin$ is not too large.

The final stage of the computation is to evaluate 
$y = \xmax + \log(1+s)$ using the computed $\shat$, for which we have
$$
  \yhat = \bigl(\xmax + \log(1+\shat)(1+\d_3)\bigr)(1+\d_4),
  \quad |\d_3|, |\d_4| \le u.
$$
Here, we are assuming that the $\logonep$ function has the property  
$$
   \fl(\logonep(s)) = \logonep(s)(1+\d), \quad |\d| \le u.
$$
Ignoring the innocuous $\d_4$ term and writing, by \eqref{werrb},
\begin{equation}\label{sb}
  \shat = s(1+\eta), \quad |\eta| \le (n + \xmax - \xmin) u + O(u^2),
\end{equation}
we have
\begin{align*}
    \yhat &= \xmax + \log(1 + s(1 + \eta))(1 + \d_3)\\
          &= \xmax + \log(1 + s + s\eta))(1 + \d_3)\\
          &= \xmax + \left(\log(1 + s) + \frac{s\eta}{1 + s} 
                           + O(u^2) \right)(1 + \d_3),
\end{align*}
using a Taylor series expansion about $1+s$ of the logarithm.
Hence
$$
   \yhat - y = \log(1+s)\d_3 + \frac{s\eta}{1 + s} (1 + \d_3) + O(u^2).
$$
Bounding $\eta$ using \eqref{sb} gives 
\begin{equation}\label{dag11}
  |y - \yhat| \le \log(1 + s)u + \frac{s}{1+s}(n+\xmax-\xmin)u +
  O(u^2)
\end{equation}
or, as a relative error bound, since $s \ge 0$,
\begin{equation}\label{dag12}
 \left|\frac{y - \yhat}{y} \right| 
  \le \left(\frac{\log(1+s) + n+\xmax-\xmin}{|y|}\right)u + O(u^2).
\end{equation}
Simplifying the bound gives the next result.

\begin{theorem}[Shifted \lse\ algorithm]\label{thm.lse}
The computed \lse\ $\yhat$ from \Alg~\ref{alg.log-sum-exp} satisfies
\begin{equation}\label{yerr}
 \left|\frac{y-\yhat}{y}\right| 
  = \left|\frac{y + n - \xmin}{y}\right|u + O(u^2). 
\end{equation}
\end{theorem}

The main question is how this result compares with Theorem~\ref{thm.lse-basic}
for the unshifted \alg.
The only difference in the bounds is that $|y|+n +1$ in 
\eqref{yerr2} is replaced by $|y+n-\xmin|$ here.
Now 
$|y+n-\xmin| \gg |y| + n$ is possible only if $\xmin \ll 0$
and $\xmin \ll \xmax$, so let us assume that these two inequalities hold.
The term $|y+n-\xmin|$ comes from bounding the 
term $(1+a-x_i) w_i$, where $w_i$ is defined in \eqref{sum}
and $x_i = \xmin$,
and if $\xmin \ll 0$ then 
$w_i = \eu^{x_i-a} = \eu^{\xmin - \xmax} \ll 1$.
Hence the potentially large constant is mitigated by the $w_i$ term
that it multiplies---something that is lost in the manipulations to achieve
a readable bound.
We conclude that shifting should have little effect on the accuracy.

We note that \eqref{yerr} is weaker than necessary when $s \ll 1$
(recall that $s \ge 0$), since we
bounded $s/(1+s)$ by $1$ in going from 
\eqref{dag11} to \eqref{dag12}. If $s \ll 1$ then \eqref{dag11}
becomes
\begin{equation*}
  |y - \yhat| \lesssim s (1 + n + \xmax - \xmin) u + O(u^2). 
\end{equation*}
Since $s\ll 1$ also implies $x_i \ll \xmax$ for $i\ne k$ and hence 
$y \approx \xmax$, we have 
\begin{equation*}
  \frac{|y - \yhat|}{|y|} \lesssim s \frac{|1 + n + y - \xmin|}{|y|} 
             u + O(u^2),
\end{equation*}
which is a factor $s$ smaller than \eqref{yerr}.



Turning to the evaluation of the \sm\ function $g$ from the shifted formula
\eqref{gshift}, we have, using \eqref{werrb},
$$
   \ghat_j = \frac{ \exp( (x_j-a) (1+\d_1) ) (1+\d_2)(1+\d_3) }{ s(1+\eta) },
$$
where $\d_2$ corresponds to the exponential evaluation and 
$\d_3$ to the division, and
$$
      \mbox{$|\d_i| \le u$, $i = 1\colon 3$}, \qquad
       |\eta| \le (n + \xmax-\xmin)u + O(u^2).
$$
Therefore
\begin{align*}
   \ghat_j &= g_j \frac{ \exp( (x_j-a) \d_1 ) (1+\d_2)(1+\d_3) }{ 1+\eta }\\
           &= g_j (1 + \th), \quad
                |\th| \le \bigl(n + 2 + 2 (\xmax-\xmin)\bigr) u + O(u^2).
\end{align*}

Hence we have obtained the following result.

\begin{theorem}[Shifted softmax algorithm]\label{thm.gshift}
The computed $\ghat$ from \Alg~\ref{alg.log-sum-exp} satisfies
\begin{equation}\label{gerr3}
  \frac{ \normo{g-\ghat} }{ \normo{g} }
  \le \bigl(n + 2 + 2(\xmax-\xmin) \bigr) u + O(u^2).
\end{equation}
\end{theorem}

Again, this is broadly commensurate with Theorem~\ref{thm.g} for the
unshifted evaluation,
bearing in mind the comments following Theorem~\ref{thm.lse}.

Finally, we consider \eqref{galt} 
with the \lse\ computed by \Alg~\ref{alg.log-sum-exp}.
In \fpa\ we have the same equation \eqref{eq:step2} as for the unshifted
\alg, but now with $\th$ bounded by, 
using \eqref{yerr},
$$
 |\th| \le (1 + |x_j-y| + |y+n-\xmin| )u + O(u^2).
$$
We have obtained the following result.

\begin{theorem}[Alternative shifted softmax algorithm]\label{thm.gshift3}
The computed $\ghat$ from \eqref{galt} with the \lse\ computed by 
\Alg~\ref{alg.log-sum-exp} satisfies
\begin{equation}\label{gerr4}
  \frac{ \normo{g-\ghat} }{ \normo{g} }
  \le \Bigl(1 + \max_j|x_j-y| + |y+n-\xmin| \Bigr)u + O(u^2).
\end{equation}
\end{theorem}

This is broadly similar to Theorem~\ref{thm.gshift2}
for the unshifted alternative \sm\ \alg.


\section{Computational experiments}\label{sec.numexp} 


We now perform some experiments in a realistic setting, 
using MATLAB R2019a.
The codes and data used for the experiments are available
online\footnote{\url{https://github.com/higham/logsumexp-softmax-tests}}.

Our aims are to examine the sharpness of the rounding error bounds and 
to give a pairwise comparison of the accuracy of the algorithms in \fpa.
Our data comes from a deep learning application.
To generate the data, we first set up and trained  
an artificial  neural network, using the MATLAB
Deep Learning Toolbox~\cite{matlab-DLT}.
More precisely, we 
trained a network to classify  
handwritten digit data
from the widely used MNIST data set 
\cite{lcb-digits}.
Here each data point is a grayscale
$28 \times 28 $ pixel image and 
there are ten categories: 
$0$, $1$, \ldots, $9$.
We used a network
whose architecture 
has the following general form:
\begin{enumerate}[nosep]
\item 
Image Input    $28\times28\times1$ with normalization.
\item 
Convolution    8 $3\times3\times1$ stride [1  1] padding 'same'.
\item 
Batch Normalization  8 channels.
\item 
ReLU
\item Max Pool $2\times2$ stride [2  2] padding [0  0  0  0].
\item 
Convolution 16 $3\times3\times8$ stride [1  1] padding 'same' .
\item 
 Batch Normalization 16 channels.
\item ReLU.
 \item 
Max Pool $2\times2$ stride [2  2] padding [0  0  0  0].
\item 
Convolution 32 $3\times3\times16$ stride [1  1] padding 'same'.
\item 
Batch Normalization 32 channels.
\item ReLU.
\item 
Fully Connected 10  layer.
\item 
Softmax.
\item 
Classification Output crossentropy.
\end{enumerate}
This is the default architecture from 
\cite{matlab-DLT}, where further details may be found.

The network was trained on 7500 images 
(750 from each of the ten categories), 
with 2500 further images
(250 from each of the ten categories) 
used for validation.

The network takes as input a $28\times 28$ matrix 
corresponding to the pixels in the image and returns a nonnegative 
$10\times 1$ vector whose $i$th component may be interpreted as the probability
that the image came from category $i$.  If we categorize according to the
highest probability from the output, then the trained network misclassifed
27 of the 2500 validation images, corresponding to a 98.9\% success rate.

The network uses single precision arithmetic, \texttt{fp32}.
In our experiments, we are concerned only with floating-point 
arithmetic issues, and we treat the 
trained network as a means to produce a realistic data set.
To do this, we extracted the 2500 single precision vectors from the
validation set that were passed into the softmax layer
and converted them to fp16 or bfloat16.
We then used this data 
in our implementation of the 
softmax and log-sum-exp algorithms 
that we have studied in the previous sections.

To record errors in computed results we applied the basic algorithm,
Algorithm~\ref{alg.basic}, in single precision 
to provide a reference solution and 
used the  \texttt{chop} function of \cite{hipr19} to simulate 
half precision arithmetic, in both the fp16 format and the bfloat16 format.

We first describe experiments in fp16.
The components in the 2500 test vectors $x \in \RR^{10}$ vary between about
$-19$ and $+20$.  As indicated in Table~\ref{table.IEEE2}, $\eu^x$ overflows
in fp16 for $x \gtrsim 11$.  Hence, in these tests, overflow is an issue
for the basic log-sum-exp implementation in Algorithm~\ref{alg.basic}: it
generated an \texttt{Inf} for 475 of the 2500 test vectors.  The shifted
version of log-sum-exp in Algorithm~\ref{alg.log-sum-exp} did not overflow.
In the plots below, we do not include results for the cases where
Algorithm~\ref{alg.basic} produced overflow.

\begin{figure}
\begin{center}
\includegraphics[scale=0.8]{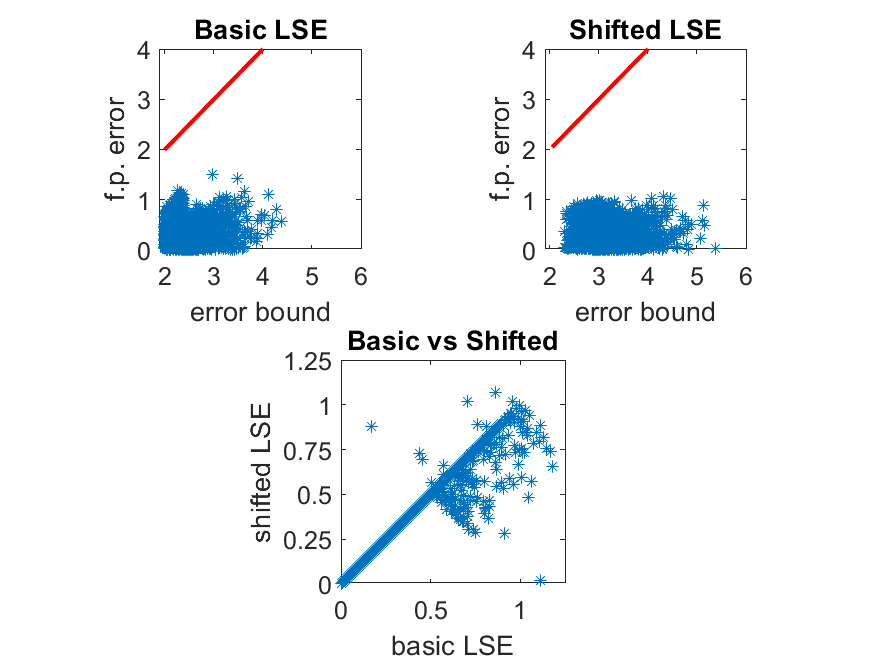}
\caption{Scatter plots of errors and error bounds, scaled by unit roundoff, 
  over $2025$ vectors in $\RR^{10}$ for \lse\ \alg s in fp16.
  See the text for a description of the axes.
  Upper left: basic implementation of
  log-sum-exp from Algorithm~\ref{alg.basic}.
  According to the error analysis, all points should
  lie below the reference line $y = x$ (shown in red).  Upper right:
  corresponding results for the shifted implementation of log-sum-exp in
  Algorithm~\ref{alg.log-sum-exp}.  
  Lower: scaled error from
  Algorithm~\ref{alg.basic} versus scaled error from
  Algorithm~\ref{alg.log-sum-exp}.}
\label{Fig:lse}
\end{center}
\end{figure}

First, we
look at the log-sum-exp algorithms.
In the upper  left plot of 
Figure~\ref{Fig:lse} 
we used the basic implementation of log-sum-exp, 
Algorithm~\ref{alg.basic}.
We scatter plot over the 2025  
vectors where no overflow occurred.
For each such vector, the horizontal coordinate is the 
leading term in the error bound of 
Theorem~\ref{thm.lse-basic}, 
scaled by $u$, that is, $1 + (n+1)/|y|$. 
Here, as shown in Table~\ref{table.IEEE},  
$u = 4.88\times 10^{-4}$ for fp16.
The vertical coordinate is the actual scaled 
relative error $| \widehat{y} - y|/(u|y|)$. 
The plot also gives a reference line of slope $1$ 
from the origin.
We see that the bound is always satisfied and is reasonably sharp in many cases.

In the upper right plot of Figure~\ref{Fig:lse} we show 
corresponding results for 
the shifted log-sum-exp implementation in 
Algorithm~\ref{alg.log-sum-exp},
using the bound from Theorem~\ref{thm.lse}.

In the lower part of Figure~\ref{Fig:lse} we scatter plot the floating-point
errors for the basic and shifted algorithms.  
Here, for 1863 out of the 2025 cases (92\%) the two errors were identical
to all digits in the half precision computation.
In more detail, over all the data points the ratio of the error in 
the basic log-sum-exp (horizontal axis) divided by the error in the 
shifted version (vertical axis) varied between 0.19 and 59, with a mean of
1.07 and a standard error of 0.03.
This indicates that the two versions perform similarly, with the shift 
producing slightly better results.    

\begin{figure}
\begin{center}
\includegraphics[scale=0.8]{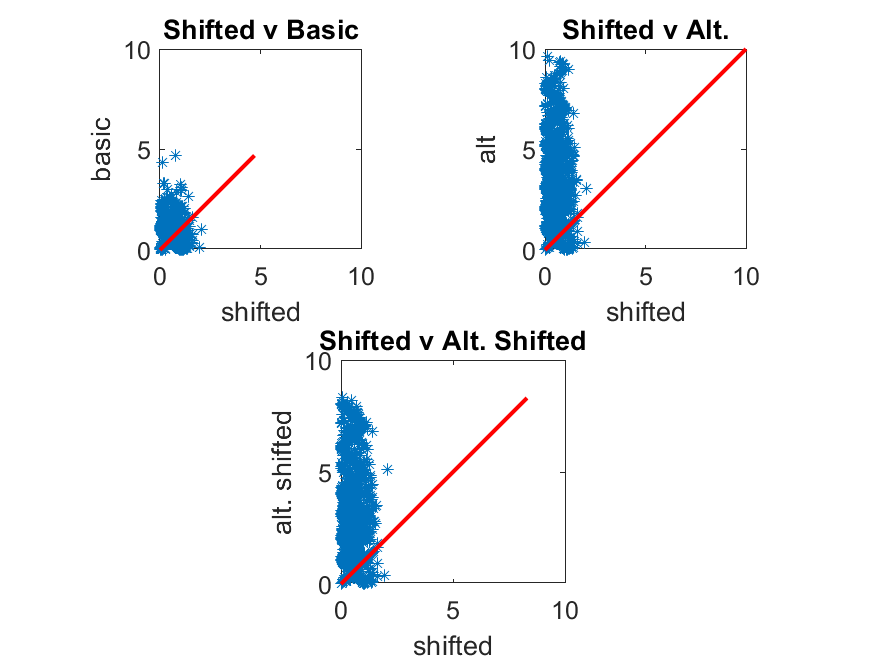}
\caption{Scatter plots of errors, scaled by unit roundoff, for
  softmax algorithms in fp16.  
  See the text for a description of the axes.
  Reference line is $y = x$.  }
\label{Fig:soft}
\end{center}
\end{figure}

We now move on to the four softmax implementations.  In
Figure~\ref{Fig:soft} we use the shifted softmax implementation from
Algorithm~\ref{alg.log-sum-exp}, analysed in Theorem~\ref{thm.gshift}, as
the basis for comparison.  The upper left plot has the scaled error 
$\normo{\ghat-g}/(u\normo{g})$
from Algorithm~\ref{alg.log-sum-exp} on the horizontal axis and the scaled error
from the basic softmax in Algorithm~\ref{alg.basic} on the vertical axis.
The upper right plot compares the shifted softmax against the alternative
algorithm analyzed in 
Theorem~\ref{thm.gshift2}. 
Similarly, the lower plot
compares against the alternative shifted softmax algorithm analyzed in
Theorem~\ref{thm.gshift3}.
We see that the \sm\ values obtained from
Algorithms~\ref{alg.basic}
and  \ref{alg.log-sum-exp} have similar accuracy,
whereas the alternative \sm\ versions based on the rewrite in (\ref{galt})
are typically less accurate.

The results in Figures~\ref{Fig:lse} and \ref{Fig:soft} are consistent with
our floating-point error analysis.

\begin{figure}
\begin{center}
\includegraphics[scale=0.55]{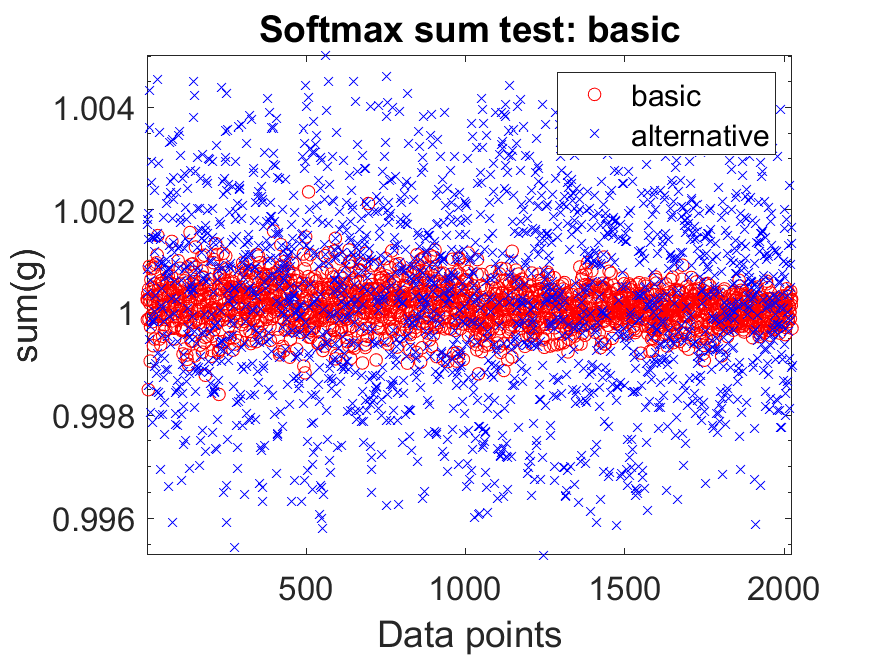}
\caption{Sum of entries of computed softmax vector for
 \Alg~\ref{alg.basic} (red circles), analyzed in Theorem~\ref{thm.g},
 and the
alternative  (blue crosses) analyzed in
  Theorem~\ref{thm.gshift2}.}
\label{Fig:sum2}
\end{center}
\end{figure}

\begin{figure}
\begin{center}
\includegraphics[scale=0.55]{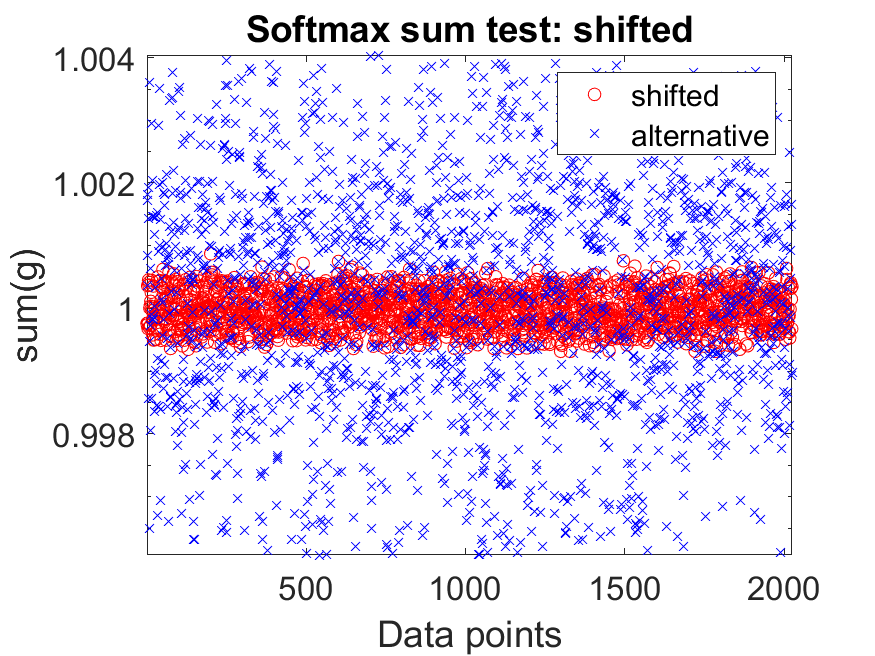}
\caption{Sum of entries of computed softmax vector for 
\Alg~\ref{alg.log-sum-exp} (red circles),
analyzed in Theorem~\ref{thm.gshift}, and the 
         alternative (blue crosses) analyzed in
         Theorem~\ref{thm.gshift3}.}
\label{Fig:sum3}
\end{center}
\end{figure}

A further test is to compute the sum of each \sm\ vector, which
should equal $1$.
In Figure~\ref{Fig:sum2} we compare the softmax sums for the basic algorithm (red circles)  
analyzed in 
Theorem~\ref{thm.g} and the alternative version (blue crosses) analyzed in
Theorem~\ref{thm.gshift2}.
Similarly, Figure~\ref{Fig:sum3} compares the shifted softmax algorithm
analyzed in Theorem~\ref{thm.gshift} and its alternative analyzed in
Theorem~\ref{thm.gshift3}.
The order along the $x$-axis is arbitrary; it corresponds to the order in
which the data vectors were generated.  These figures provide further
evidence that the alternative \sm\ \alg s are less accurate
than the basic or shifted \alg s.

We also conducted the corresponding experiments in simulated bfloat16
arithmetic.  Here, as indicated in Tables~\ref{table.IEEE} and
\ref{table.IEEE2}, the number range is increased at the expense of reduced
precision.  In this case there was no overflow in any of the algorithms.
The results were very similar to those for fp16, so they are not shown~here.





\section{Conclusions}\label{sec.conc}

The \lse\ and \sm\ functions both feature in many computational
pipelines, so it is important to compute them accurately and to avoid 
generating infs or NaNs because of overflow or underflow.
To this end, a shift is usually incorporated into the defining formulas,
yielding \eqref{ynew} and \eqref{gshift}.  It is important to understand
the effect of the shift on the accuracy of the computed result, especially
when computations are carried out in a low precision such as bfloat16 or
fp16, which have the equivalent of only 3 or 4 decimal digits of precision.

Our \rea\ shows that shifting by the largest element of the input vector
does not lessen the accuracy of the computed
\lse\ and \sm.  
Underlying this pleasing fact is the phenomenon that 
any large coefficients caused by shifting are canceled by 
multiplication with small exponentials.

We obtained an explicit formula for the \cn\ of \lse\ and 
bounds for the \cn\ of \sm, 
and we were able to identify situations in which the \lse\ \alg s are
guaranteed to be forward stable.

For the alternative and widely used \sm\ formula that avoids division,
\eqref{galt}, we obtained larger error bounds than for the shifted
formula~\eqref{gshift}.  Since our numerical experiments confirm that
larger errors are typically obtained in practice, we recommend using
\eqref{gshift} instead of~\eqref{galt} to evaluate \sm.

In summary, \Alg~\ref{alg.log-sum-exp} is our recommendation for computing
\lse\ and softmax.  It avoids overflow, reduces the chance of 
\new{harmful} underflow,
and generally produces results as accurate as those from the unshifted formulas.

\bibliographystyle{plain}
\bibliography{strings,paper,njhigham}

\def\noopsort#1{}\def\hbk{hardback}\def\pbk{paperback}
\begin{thebibliography}{10}

\bibitem{aphi14}
Mary Aprahamian and Nicholas~J. Higham.
\newblock \href{https://doi.org/10.1137/130920137}{The matrix unwinding
  function, with an application to computing the matrix exponential}.
\newblock {\em SIAM J. Matrix Anal. Appl.}, 35\penalty0 (1):\penalty0 88--109,
  2014.

\bibitem{bish06}
Christopher~M. Bishop.
\newblock {\em Pattern Recognition and Machine Learning}.
\newblock Spring{\-}er-Ver{\-}lag, New York, 2006.
\newblock xx+738 pp.
\newblock ISBN 978-0-387-31073-2.

\bibitem{bova04}
Stephen Boyd and Lieven Vandenberghe.
\newblock \href{http://web.stanford.edu/~boyd/cvxbook}{{\em Convex
  Optimization}}.
\newblock Cambridge University Press, Cambridge, UK, 2004.
\newblock xiii+716 pp.
\newblock ISBN 0-521-83378-7.

\bibitem{cgp19}
Giuseppe~C. Calafiore, Stephane Gaubert, and Corrado Possieri.
\newblock \href{https://doi.org/10.1109/tnnls.2019.2910417}{Log-sum-exp neural
  networks and posynomial models for convex and log-log-convex data}.
\newblock {\em IEEE Trans. Neural Networks and Learning Systems}, pages 1--12,
  2019.

\bibitem{dlm07}
Florent de~Dinechin, Christoph Lauter, and Jean-Michel Muller.
\newblock \href{https://doi.org/10.1051/ita:2007003}{Fast and correctly rounded
  logarithms in double-precision}.
\newblock {\em RAIRO-Inf. Theor. Appl.}, 41\penalty0 (1):\penalty0 85--102,
  2007.

\bibitem{matlab-DLT}
{Deep Learning Toolbox}.
\newblock The MathWorks, Inc., Natick, MA, USA.
\newblock \url{http://www.mathworks.co.uk/products/deep-learning/}.

\bibitem{efha16}
Bradley Efron and Trevor Hastie.
\newblock \href{http://dx.doi.org/10.1017/CBO9781316576533}{{\em Computer Age
  Statistical Inference. {Algorithms}, Evidence, and Data Science}}.
\newblock Cambridge University Press, Cambridge, UK, 2016.
\newblock xix+475 pp.
\newblock ISBN 978-1-107-14989-2.

\bibitem{gbc16}
Ian Goodfellow, Yoshua Bengio, and Aaaron Courville.
\newblock {\em Deep Learning}.
\newblock The MIT Press, Cambridge, MA, USA, 2016.
\newblock xxii+775 pp.
\newblock ISBN 978-0-262-03561-3.

\bibitem{hp82}
{\em {HP-\textup{15}C} Advanced Functions Handbook}.
\newblock Hewlett-Packard, Portable Computer Division, Corvallis, OR, USA,
  1982.
\newblock 221 pp.
\newblock Part number 00015-90011 Rev.~C.

\bibitem{hh19}
Catherine~F. Higham and Desmond~J. Higham.
\newblock Deep learning: {A}n introduction for applied mathematicians.
\newblock {\em SIAM Review}, to appear, 2019.

\bibitem{high93s}
Nicholas~J. Higham.
\newblock \href{https://doi.org/10.1137/0914050}{The accuracy of floating point
  summation}.
\newblock {\em SIAM J. Sci. Comput.}, 14\penalty0 (4):\penalty0 783--799, 1993.

\bibitem{high:ASNA2}
Nicholas~J. Higham.
\newblock \href{http://dx.doi.org/10.1137/1.9780898718027}{{\em Accuracy and
  Stability of Numerical Algorithms}}.
\newblock Second edition, Society for Industrial and Applied Mathematics,
  Philadelphia, PA, USA, 2002.
\newblock xxx+680 pp.
\newblock ISBN 0-89871-521-0.

\bibitem{high:FM}
Nicholas~J. Higham.
\newblock \href{http://dx.doi.org/10.1137/1.9780898717778}{{\em Functions of
  Matrices: {Theory} and Computation}}.
\newblock Society for Industrial and Applied Mathematics, Philadelphia, PA,
  USA, 2008.
\newblock xx+425 pp.
\newblock ISBN 978-0-898716-46-7.

\bibitem{hima18a}
Nicholas~J. Higham and Theo Mary.
\newblock \href{http://eprints.maths.manchester.ac.uk/2689/}{A new approach to
  probabilistic rounding error analysis}.
\newblock {MIMS EPrint} 2018.33, Manchester Institute for Mathematical
  Sciences, The University of Manchester, UK, November 2018.
\newblock 22 pp.
\newblock Revised March 2019. To appear in SIAM J. Sci. Comput.

\bibitem{hipr19}
Nicholas~J. Higham and Srikara Pranesh.
\newblock \href{http://eprints.maths.manchester.ac.uk/2723/}{Simulating low
  precision floating-point arithmetic}.
\newblock {MIMS EPrint} 2019.4, Manchester Institute for Mathematical Sciences,
  The University of Manchester, UK, March 2019.
\newblock 18 pp.
\newblock Revised July 2019. To appear in SIAM J. Sci. Comput.

\bibitem{ieee08}
\href{http://dx.doi.org/10.1109/IEEESTD.2008.4610935}{{\em {IEEE} Standard for
  Floating-Point Arithmetic, {IEEE} {Std} 754-2008 (revision of {IEEE} Std
  754-1985)}}.
\newblock IEEE Computer Society, New York, 2008.
\newblock 58 pp.
\newblock ISBN 978-0-7381-5752-8.

\bibitem{inte18}
{Intel Corporation}.
\newblock
  \href{https://software.intel.com/en-us/download/bfloat16-hardware-numerics-definition}{{BFLOAT16}---hardware
  numerics definition}, November 2018.
\newblock White paper. Document number 338302-001US.

\bibitem{scipy}
Eric Jones, Travis Oliphant, Pearu Peterson, et~al.
\newblock {SciPy}: Open source scientific tools for {Python}, 2001--.
\newblock \url{http://www.scipy.org/}.

\bibitem{lcb-digits}
Yann LeCun, Corinna Cortes, and Christopher J.~C. Burges.
\newblock \href{http://yann.lecun.com/exdb/mnist/}{The {MNIST} database of
  handwritten digits}.
\newblock Accessed June 17, 2019.

\bibitem{matlab-prml}
Matlab code for machine learning algorithms in book {PRML}.
\newblock \url{https://github.com/PRML/PRMLT}.

\bibitem{mull16}
Jean-Michel Muller.
\newblock \href{http://dx.doi.org/10.1007/978-1-4899-7983-4}{{\em Elementary
  Functions: {Algorithms} and Implementation}}.
\newblock Third edition, Birkh{\"{a}}user, Boston, MA, USA, 2016.
\newblock xxv+283 pp.
\newblock ISBN 978-1-4899-7981-0.

\bibitem{mbdj18}
Jean-Michel Muller, Nicolas Brunie, Florent de~Dinechin, Claude-Pierre
  Jeannerod, Mioara Joldes, Vincent Lef{\`e}vre, Guillaume Melquiond, Nathalie
  Revol, and Serge Torres.
\newblock \href{http://dx.doi.org/10.1007/978-3-319-76526-6}{{\em Handbook of
  Floating-Point Arithmetic}}.
\newblock Second edition, Birkh{\"{a}}user, Boston, MA, USA, 2018.
\newblock xxv+627 pp.
\newblock ISBN 978-3-319-76525-9.

\bibitem{murp12}
Kevin~P. Murphy.
\newblock {\em Machine Learning: a Probabilistic Approach}.
\newblock Cambridge University Press, Cambridge, UK, 2012.

\bibitem{matlab-SMLT}
{Statistics and Machine Learning Toolbox}.
\newblock The MathWorks, Inc., Natick, MA, USA.
\newblock \url{https://uk.mathworks.com/products/statistics.html}.

\bibitem{rlang}
{R}~Core Team.
\newblock \href{https://www.R-project.org/}{{\em {R}: {A} Language and
  Environment for Statistical Computing}}.
\newblock R Foundation for Statistical Computing, Vienna, Austria.

\bibitem{wiba98}
Christopher K.~I. Williams and David Barber.
\newblock \href{https://doi.org/10.1109/34.735807}{Bayesian classification with
  {Gaussian} processes}.
\newblock {\em IEEE Trans. Pattern Analysis and Machine Intelligence},
  20\penalty0 (12):\penalty0 1342--1351, 1998.

\end{thebibliography}

\end{document}